\documentclass[10pt,twoside]{article}
\usepackage{graphicx}
\usepackage{amsmath}
\usepackage{Latex-document}

\title{\bf  The Dynamical Systems Approach to the\vskip -2mm Equations
of a Linearly Viscous\vskip -2mm Compressible Barotropic Fluid
\vskip 6 mm}

\author{E. Feireisl\vspace*{-0.5cm}\thanks{Institute of Mathematics,
Czech Academy of Sciences, \v Zitn\' a 25, 115 67 Praha 1,
Czech Republic. E-mail: feireisl@math.cas.cz}}

\date{\vspace{-8mm}}

\begin{document}

\maketitle \thispagestyle{first} \setcounter{page}{295}

\begin{abstract}

\vskip 3mm

We develop a dynamical systems theory for the compressible Navier-Stokes equations based on global in time weak
solutions. The following questions will be addressed:

\begin{itemize}
\item
Global existence and critical values of the adiabatic constant;

\item
dissipativity in the sense of Levinson - bounded absorbing sets;

\item
asymptotic compactness;

\item
the long-time behaviour and attractors.
\end{itemize}
\vskip 4.5mm


\end{abstract}

\vskip 12 mm

\section{Introduction} \label{i}\setzero
\vskip-5mm \hspace{5mm}

The long-time behaviour of solutions to the evolutionary equations
arising in the mathematical fluid mechanics has been the subject
of many theoretical studies. This type of problems is apparently
related to the phenomena of turbulence, and there is still a
significant gap between many formal ``scenarios" and
mathematically rigorous results.

The dynamical systems in question are usually related to a system
of partial differential equations and, consequently, they are
defined in an infinite dimensional phase space. On the other
hand, the presence of dissipative terms in the equations due to
viscosity results in the existence of global attractors-compact
invariant sets attracting uniformly in time all trajectories
emanating from a given bounded set. Such a theory is well
developed for the incompressible linearly viscous fluids, and the
reader may consult the monographs of BABIN and VISHIK
\cite{BAVI}, TEMAM \cite{TE} or CONSTANTIN et al. \cite{CFT} for
the recent state of art.

On the other hand, much less seems to be known for the
compressible fluids. While there is an existence theory of the
weak solutions for the incompressible Navier-Stokes equations due
to LERAY \cite{LER}, its ``compressible" counterpart appeared
only recently in the work of LIONS \cite{LI4}. Even in the
incompressible case, there is a qualitative difference between
the two-dimensional case solved by LADYZHENSKAYA \cite{LAD}, and
the three-dimensional case representing one of the most
challenging unsolved problems of modern mathematics. It is
worth-noting that a similar gap divides the one and
more-dimensional problems for the compressible fluids.

The time evolution of the fluid density $\varrho = \varrho (t,x)$
and the velocity $\vec{u} = \vec{u}(t,x)$ is governed by the
Navier-Stokes system of equations:

\begin{equation} \label{i1}
\partial_t \varrho + {\rm div}(\varrho \vec{u}) = 0,
\end{equation}

\begin{equation} \label{i2}
\partial_t (\varrho \vec{u}) + {\rm div}(\varrho \vec{u} \otimes \vec{u} ) +
\nabla p = {\rm div} \ {S} + \varrho \vec{f},
\end{equation}
where $p$ is the pressure, ${S}$ the viscous stress tensor, and
$\vec{f}$ a given external force.

In what follows, we consider linearly viscous (Newtonian) fluids
where the viscous stress is related to the velocity by the
constitutive law
\begin{equation} \label{i3}
{S} = {\mu} \Big( \nabla \vec{u} + \nabla \vec{u}^T \Big)
+ \lambda \ {\rm div} \ \vec{u} \ {I},
\end{equation}
where the viscosity coefficients satisfy
\begin{equation} \label{i4}
\mu > 0 , \ \lambda + \mu \geq 0.
\end{equation}

Generally speaking, the pressure $p$ depends on the density and the
internal energy (temperature) of the fluid. If it is the case, the
system (\ref{i1}), (\ref{i2}) is not closed and should be
complemented by the energy equation. Unfortunately, however, the
available global existence results for this full system allow for
only for small initial data (cf. MATSUMURA and NISHIDA \cite{MANI}).

There are physically relevant situations when one can assume the
flow is barotropic, i.e., the pressure depends solely on the density.
This is the case when either the temperature (the isothermal case)
or the entropy (the isentropic case) are supposed to be constant.
The typical constitutive relation between the pressure and the
density then reads
\begin{equation} \label{i5}
p = p(\varrho) = a \varrho^{\gamma} , \ a > 0,
\end{equation}
where $\gamma = 1$ in the isothermal case, and $\gamma > 1$
represents the adiabatic constant in the isentropic regime.
More general and even non-monotone pressure-density constitutive laws
arise in nuclear plasma physics (see \cite{DFPS}). In the barotropic
regime, the equations (\ref{i1}), (\ref{i2}) form a closed system and
complemented by suitable initial and boundary conditions represent
a (at least formally) well-posed problem.

If the problem is posed on a spatial domain $\Omega \subset R^N$,
one usually assumes that the fluid adheres completely to the
boundary which is mathematically expressed by the no-slip boundary
conditions for the velocity:
\begin{equation} \label{i6}
\vec{u}|_{\partial \Omega} = 0.
\end{equation}
Note that for viscous fluids such a condition is in a very good
agreement with physical experiment.

In accordance with the deterministic principle, the state of the
system at any time $t > t_0$ should be given by the initial conditions
\begin{equation} \label{i7}
\varrho(t_0) = \varrho_I , \ (\varrho \vec{u})(t_0) = \vec{q}_I .
\end{equation}
The function $\varrho_I$ is non-negative and the momentum
$\vec{q}_I$ satisfies the compatibility condition
\begin{equation} \label{i8}
\vec{q}_I = 0  \ \mbox{a.a. on the set}\ \{ \varrho_I = 0 \} .
\end{equation}
The reason why we impose the initial conditions for the
momentum $\varrho \vec{u}$ rather than for the velocity $\vec{u}$ is
that the former quantity is weakly continuous with respect to time
while the instantaneous values of the velocity are determined only
almost anywhere with respect to time. Clearly such a problem does not
arise provided the initial density $\varrho_I$ is strictly positive.
However, it is an interesting open problem whether or not this
property is preserved at any positive time for any distributional
solution of the problem satisfying the natural energy estimates
(cf. HOFF and SMOLLER \cite{HOSM}).

\section{Finite energy weak solutions and well-posedness} \label{f}
\setzero
\vskip-5mm \hspace{5mm }

In order to study the long-time behaviour, one should first make sure
that the class of objects one deals with is not void. More precisely,
one should be able to prove the existence of global-in-time solution
for any initial data $\varrho_I$, $\vec{q}_I$ satisfying some
physically relevant hypothesis.

Multiplying the equations of motion by $\vec{u}$ and integrating by
parts one deduces the energy inequality
\begin{equation} \label{f1}
{{\rm d} \over {\rm d}t} E[\varrho, \vec{u}] +
\int_{\Omega} \mu |\nabla \vec{u} |^2 + (\lambda + \mu) |{\rm div} \ \vec{u}
|^2 \ {\rm d}x \leq \int_{\Omega} \varrho \vec{f} \cdot \vec{u} \ {\rm d}x ,
\end{equation}
where the total energy $E$ is given by the formula
\[
E[\varrho, \vec{u} ] = \int_{\Omega} \varrho |\vec{u}|^2 +
P(\varrho) \ {\rm d}x
\]
with
\[
P'(\varrho) \varrho - P(\varrho) = p(\varrho).
\]

Note that in the most common isentropic case, the function $P$ can be
taken in the form
\[
P(\varrho) = {a \over \gamma - 1} \varrho^{\gamma} .
\]

The energy inequality is the main (and almost the only one)
source of a priori estimates. Accordingly, ``reasonble" solutions
of the problem (\ref{i1}), (\ref{i2}) defined on a bounded time
interval $I \subset R$ should belong to the class
\begin{equation} \label{f2}
\varrho \geq 0, \ \varrho \in L^{\infty}(I; L^{\gamma}(\Omega)),\
\vec{u} \in L^2(I; W^{1,2}_0 (\Omega, R^N)).
\end{equation}

The energy inequality (\ref{f1}) represents an additional constraint
imposed on any solution $\varrho$, $\vec{u}$ of the problem. Similarly
as in the theory of the variational (weak) solutions developed for
the incompressible case by Leray, it is not clear if it is satisfied
for any weak solution of the problem.

Following DiPERNA and LIONS \cite{DL} we shall say that $\varrho$,
$\vec{u}$ is a renormalized solution of the continuity equation
(\ref{i1}) if the identity
\begin{equation} \label{f3}
\partial_t b(\varrho) + {\rm div} (b(\varrho) \vec{u}) +
\Big( b'(\varrho) \varrho - b(\varrho) \Big) {\rm div} \ \vec{u} = 0
\end{equation}
holds in the sense of distributions for any function
$b \in C^1(R)$ such that
\begin{equation} \label{f4}
b'(\varrho) = 0 \ \mbox{for all} \ \varrho \geq {\rm {const}}(b) .
\end{equation}
Similarly as for the energy inequality, it is not known if any weak
solution $\varrho$, $\vec{u}$ of (\ref{i1}) satisfies (\ref{f3}).

{\bf Definition 2.1} {\it
We shall say that $\varrho$, $\vec{u}$ is
a finite energy weak solution of the problem (\ref{i1} - \ref{i6})
on a set $I \times \Omega$ if the following conditions are satisfied:
\begin{itemize}
\item
The functions $\varrho$, $\vec{u}$ belong to the function spaces
determined in (\ref{f2});
\item
the energy inequality (\ref{f1}) is satified in ${\cal D}'(I)$;
\item
the continuity equation (\ref{i1}) as well as its renormalized
version (\ref{f3}) hold in ${\cal D}'(I \times R^N)$ provided $\varrho$,
$\vec{u}$ were extended to be zero outside $\Omega$;
\item
the momentum equation (\ref{i2}) is satisfied in ${\cal D}'(I \times
\Omega)$.
\end{itemize}
}

The most general available existence result reads as follows:

{\bf Theorem 2.1} {\it Let $\Omega \subset R^N$, $N=2,3$ be a bounded Lipschitz domain. Let $I = (0,T)$, and let
the initial data $\varrho_I$, $\vec{q}_I$ satisfy (\ref{i8}) together with
\[
\varrho_I \in L^{\gamma}(\Omega),
\ {|\vec{q}_I |^2 \over \varrho_I} \in L^1(\Omega).
\]
Let $\vec{f}$ be a bounded measurable function of $t \in I$,
$x \in \Omega$. Finally, let the pressure $p \in C[0,\infty)
\cap C^1(0,\infty)$ be given by a
constitutive law
\[
p = p(\varrho) , \ {1 \over a} \varrho^{\gamma - 1} - b
\leq p'(\varrho) \leq a \varrho^{\gamma - 1} + b,\
\mbox{for all}\ \varrho > 0,
\]
where $a > 0$, $b \geq 0$, and
\[
\gamma > {N \over 2}.
\]

Then the problem (\ref{i1} - \ref{i6}) admits a finite energy weak
solution $\varrho$, $\vec{u}$ on $I \times \Omega$ satisfying the
initial conditions (\ref{i7}).
}

In his pioneering work, LIONS \cite{LI4} proved Theorem 2.1
for $\Omega$ regular, $p$ monotone, and $\gamma \geq {3 \over 2}$
if $N=2$, and $\gamma \geq {9 \over 5}$ for $N=3$. The hypotheses
concerning $\gamma$ were relaxed in \cite{EF54}, \cite{FNP}, the case of a
general bounded domain $\Omega$ treated in \cite{FNP1}, and the
hypothesis of monotonicity of the pressure removed in
\cite{EF61}, \cite{DFPS}.

\section{Ultimate boundedness} \label{u}
\setzero
\vskip-5mm \hspace{5mm }

The first issue to be discussed when describing the long-time
asymptotics of a given dynamical system is ultimate boundedness or
dissipativity. This means there exists an absorbing set bounded
in a suitable topology. Here ``suitable topology" is of course
that one induced by the total energy $E$. One of possible results
in this direction is contained in the following theorem.

{\bf Theorem 3.1} {\it
Let $\Omega \subset R^N$, $N=2,3$ be a bounded
Lipschitz domain. Let $\bf{f}$ be a bounded measurable function such
that
\[
{\rm ess} \sup_{t \in R, x \in \Omega} | \vec{f}(t,x) | \leq F.
\]
Let the pressure $p$ be given by the isentropic constitutive law
\[
p= a \varrho^{\gamma} \ \mbox{with}\ \gamma > 1
\ \mbox{if}\ N=2, \ \gamma > {5 \over 3} \ \mbox{for}\ N=3.
\]
Finally, set
\[
\int_{\Omega} \varrho \ {\rm d}x = m > 0.
\]

Then there exists a constant $E_{\infty}$ depending solely on
$m$ and $F$ having the following property:

Given $E_I$, there exists a time $T=T(E_I)$ such that
\[
E[\varrho, \vec{u}] (t) \leq E_{\infty} \ \mbox{for a.a.}
t > T
\]
provided
\[
{\rm ess} \limsup_{t \to 0+} E[\varrho, \vec{u}](t) \leq E_I
\]
and $\varrho$, $\vec{u}$ is a finite energy weak solution of the
problem (\ref{i1} - \ref{i6}) on $(0,\infty) \times \Omega$.
}

The proof of Theorem 3.1 can be found in \cite{FP14} and
\cite{EF57}. The reader will have noticed the ``critical" exponent
$\gamma > {5 \over 3}$ which is larger than in the existence
Theorem 2.1 and, as a matter of fact, does not include any
physically relevant case. Indeed the value ${5 \over 3}$ happens
to be the largest adiabatic constant corresponding to a
monoatomic gas.

Under the hypotheses of Theorem 3.1, the dissipative mechanism
induced by viscocity is strong enough to prevent any ``resonance"
phenomena, i.e., the existence of unbounded solutions driven by a
bounded external force. Another interesting feature is the
existence of periodic solutions provided $\vec{f}$ is periodic in
time.

{\bf Theorem 3.2} {\it
In addition to the hypotheses of Theorem 3.1, assume that
$\vec{f}$ is time periodic, i.e.,
\[
\vec{f}(t + \omega, x) = \vec{f}(t,x) \ \mbox{for a.a.}\
t \in R , \ x \in \Omega
\]
with a certain period $\omega > 0$.

Then there exists at least one finite energy weak solution of the
problem (\ref{i1} - \ref{i6}) on $R \times \Omega$ which is
$\omega-$periodic in time, i.e.,
\[
\varrho(t + \omega) = \varrho(t),\ (\varrho \vec{u})(t + \omega) =
(\varrho \vec{u})(t) \ \mbox{for all}\ t \in R.
\]
}

See \cite{FMPS} for the proof.

\section{Asymptotic compactness} \label{a}
\setzero
\vskip-5mm \hspace{5mm }

The property of asymptotic compactness of a given dynamical system
plays a crucial role in the proof of existence of a global (compact)
attractor. Here, the main problem is the density component which is
bounded only in $L^{\gamma}(\Omega)$, and, consequently, compact only
with respect to the weak topology on this space. Moreover, given the
hyperbolic character of the continuity equation, one cannot hope
any possible oscillations of the density to be killed at a finite
time. However, the amplitude of possible oscillations is decreasing
in time uniformly on trajectories emanating from a given bounded set.

{\bf Theorem 4.1} {\it
Let $\Omega \subset R^N$, $N=2,3$ be a bounded Lipschitz domain.
Let the pressure $p$ be given by the isentropic constitutive relation
\[
p = a \varrho^{\gamma}, \ a > 0, \ \gamma > 1 \ \mbox{if}\
N=2, \ \gamma > {5\over 3} \ \mbox{for}\ N=3.
\]
Let $\vec{f}_n$ be a sequence of functions such that
\[
{\rm ess} \sup_{t \in R, \ x \in \Omega} |\vec{f}_n(t,x)| \leq F
\ \mbox{independently of}\ n=1,2,...\ .
\]
Finally, let $\varrho_n$, $\vec{u}_n$ be a sequence of finite energy
weak solutions to the problem (\ref{i1} - \ref{i6}) on $(0,\infty)
\times \Omega$ such that
\[
\int_{\Omega} \varrho_n \ {\rm d}x = m,
\]
\[
{\rm ess} \limsup_{t  \to 0+} E[\varrho_n , \vec{u}_n ] (t) \leq E_I
\]
independently of $n=1,2,...\ .$

Then any sequence of times $t_n \to \infty$ contains a subsequence
such that
\[
\varrho_n (t_n + t) \to \varrho (t) \ \mbox{strongly in}
\ L^1(\Omega) \ \mbox{and weakly in}\ L^{\gamma}(\Omega),
\]
\[
(\varrho_n \vec{u}_n) (t_n + t) \to (\varrho \vec{u})(t)
\ \mbox{weakly in}\ L^1(\Omega,R^N)
\]
for any $t \in R$.

Moreover,
\[
\int_J \int_{\Omega} | \varrho_n (t_n + t,x) - \varrho(t,x)
|^{\gamma} \ {\rm d}xdt \to 0,
\]
\[
\int_{J} \int_{\Omega} | (\varrho_n \vec{u}_n)(t_n + t,x) -
(\varrho \vec{u})(t,x) | \ {\rm d}xdt \to 0
\]
for any bounded interval $J \subset R$. The limit functions
$\varrho$, $\vec{u}$ represent a globally defined (for $t \in R$)
finite energy weak solution of the problem (\ref{i1} - \ref{i6})
driven by a force
\[
\vec{f} = \lim_{n \to \infty} \vec{f}_n (t_n + \cdot)
\ \mbox{in the weak star topology of the space}\
L^{\infty} (R \times \Omega,R^N),
\]
and such that
\[
{\rm ess} \sup_{t \in R} E[\varrho, \vec{u}] (t) < \infty .
\]
}

The proof of Theorem 4.1 is given in \cite{FP15} (see also \cite{EF53}).

\section{The long-time behaviour, attractors} \label{l}
\setzero
\vskip-5mm \hspace{5mm }

The results presented in Sections \ref{u}, \ref{a} allow us to develop a theory of attractors analogous to that
one for the incompressible flows (see e.g. TEMAM \cite{TE}). Assume, for the sake of simplicity, that the driving
force $\vec{f}$ is independent of $t$. We introduce
$$
{\cal A} = \{ [\varrho_I, \vec{q}_I ] \ | \ \varrho_I = \varrho(0),\ \vec{q}_I = (\varrho \vec{u})(0) \
\mbox{where} \ \varrho,\ \vec{u} \ \mbox{is a finite energy weak solution}
$$
\begin{equation} \label{l1}
\mbox{of the problem (\ref{i1} - \ref{i6}) on}\ R \times \Omega \ \mbox{with}\ E[\varrho, \vec{u}] \in
L^{\infty}(R) \}.
\end{equation}

In other words, the set ${\cal A}$ is formed by all globally
defined (for $t \in R$) trajectories whose energy is uniformly
bounded.

The next result shows that ${\cal A}$ is a global attractor in the
sense of FOIAS and TEMAM \cite{FOTE}.

{\bf Theorem 5.1} {\it
Let $\Omega \subset R^N$, $N=2,3$ be a bounded Lipschitz domain.
Let $\vec{f} = \vec{f}(x)$ be a bounded measurable
function independent of time, and let the pressure
$p$ be given by the isentropic constitutive law
\[
p = a \varrho^{\gamma},\ a > 0, \ \gamma > 1
\ \mbox{for}\ N=2,\ \gamma > {5\over 3} \ \mbox{if}\
N=3.
\]

Then the set ${\cal A}$ defined by (\ref{l1}) is compact in the space
\[
L^{\alpha}(\Omega) \times L^{2 \gamma \over \gamma + 1}_{weak}
(\Omega),
\]
and
\[
\sup_{[\varrho ,  \vec{u}] \in {\cal B}(E_I) }
\Big[ \inf_{[\varrho_I, \vec{q}_I] \in {\cal A}}
\Big( \| \varrho(t) - \varrho_I \|_{L^{\alpha}(\Omega)} +
\Big| \int_{\Omega} ((\varrho \vec{u})(t) - \vec{q}_I ) \cdot \phi
\ {\rm d}x \Big| \Big) \Big] \to 0
\]
\[
\ \mbox{for}\ t \to \infty
\]
for any $1 \leq \alpha < \gamma$ and any
$\phi \in
L^{2 \gamma \over \gamma - 1}(\Omega,R^N)$, where
the symbol ${\cal B}(E_I)$ stands for the set of all finite energy
weak solutions of the problem (\ref{i1} - \ref{i6}) on $(0,\infty)
\times \Omega$ such that
\[
{\rm ess} \limsup_{t \to 0} E[\varrho, \vec{u}] (t) \leq E_I.
\]

}

The proof can be found in \cite{EF56}.

To conclude, we give another result concerning the long-time
behaviour of solutions on condition that the driving force is a
gradient of a scalar potential.

{\bf Theorem 5.2} {\it
Let $\Omega \subset R^N$, $N=2,3$ be a bounded Lipschitz domain.
Let the pressure $p$ be given by the isentropic constitutive relation
\[
p = a \varrho^{\gamma},\ a > 0, \ \gamma > {N \over 2}.
\]
Let the driving force $\vec{f}$ be of the form
\[
\vec{f} = \nabla F,
\]
where $F = F(x)$ be a scalar potential which is globally Lipschitz on
$\overline{\Omega}$ and such that the upper level sets
\[
[F > k] = \{ x \in \Omega \ | \ F(x) > k \}
\]
are connected for any $k \in R$.

Then any finite energy weak solution $\varrho$, $\vec{u}$ of the
problem (\ref{i1} - \ref{i6}) satisfies
\[
\varrho (t) \to \varrho_s \ \mbox{in}\ L^{\gamma}(\Omega),
\ (\varrho \vec{u})(t) \to 0 \ \mbox{in}\
L^1(\Omega) \ \mbox{as}\ t \to \infty ,
\]
where $\varrho_s$ is a solution of the static problem
\[
a \nabla \varrho_s^{\gamma} = \varrho_s \vec{f} \ \mbox{on}
\ \Omega .
\]
}

{\bf Acknowledgement } The work was supported by Grant No.
201/02/0854 of GA \v CR.

\label{lastpage}

\end{document}